\documentclass[12 pt]{amsproc}
\usepackage{geometry} 
\usepackage{fullpage}
\usepackage{graphicx}
\usepackage{amssymb}
\usepackage{epstopdf}
\usepackage{amsmath,amsthm,amscd,amssymb}
\usepackage{latexsym}
\usepackage[colorlinks,citecolor=red,pagebackref,hypertexnames=false]{hyperref}

\numberwithin{equation}{section}

\theoremstyle{plain}
\newtheorem{theorem}{Theorem}[section]
\newtheorem{lemma}[theorem]{Lemma}
\newtheorem{corollary}[theorem]{Corollary}

\theoremstyle{definition}

\theoremstyle{remark}
\newtheorem{remark}[theorem]{Remark}

\newtheorem{case[theorem]}{Case}

\def \R{{\mathbb R}}

\def\norm#1.#2.{\lVert#1\rVert_{#2}}

\def\R{\mathbb R}

\def\hdim{\hbox{dim}_{\mathcal H}}

\title{An elementary 
approach to simplexes \\ in thin subsets of Euclidean space}   
\date{August 17, 2016}             
\author{A. Greenleaf, A. Iosevich, B. Liu and E. Palsson}

\email{allan@math.rochester.edu} 
\email{iosevich@math.rochester.edu}
\email{bliu19@ur.rochester.edu}
\email{palsson@vt.edu}

\address{Department of Mathematics, University of Rochester, Rochester, NY 14627}

\address{Department of Mathematics,
Viriginia Tech University, Blacksburg, VA 24061}

\thanks{The  first listed author was partially supported by NSF Grant DMS-1362271 and a Simon Foundation Fellowship,
and the second listed author by NSA Grant H98230-15-1-0319.}

\begin{document}
\maketitle

\begin{abstract} We prove that if the Hausdorff dimension of $E \subset {\Bbb R}^d$, $d \ge 3$, is greater than
$\min \left\{ \frac{dk+1}{k+1}, \frac{d+k}{2} \right\},$ then the ${k+1 \choose 2}$-dimensional Lebesgue measure of $T_k(E)$, the set of congruence classes of $k$-dimensional simplexes with vertices in $E$, is positive. This improves the best  bounds previously known, decreasing the $\frac{d+k+1}{2}$ threshold obtained in \cite{EHI12} to $\frac{d+k}{2}$ via a  different and conceptually simpler method. 
We also give a simpler  proof of the $d-\frac{d-1}{2d}$ threshold for $d$-dimensional simplexes obtained in \cite{GI12,GGIP12}. 

\end{abstract}

\section{Introduction} 
\label{introduction} 

The classical Falconer distance problem, introduced by Falconer in \cite{Fal86} (see also \cite{Mat95} for the background information) is to find 
the dimensional threshold $s_0=s_0(d)$ such that if the Hausdorff dimension of a compact set $E \subset {\Bbb R}^d$, $d \ge 2$, is greater 
than $s_0$, then the Lebesgue measure of $\Delta(E):=\{|x-y|: x,y \in E \}$ is positive. Here, and throughout, $|x|=\sqrt{x_1^2+\dots+x_d^2}$, 
the usual Euclidean distance. The best  result known in the general setting, due to Wolff \cite{W99} for $d=2$ and  Erdo\u{g}an \cite{Erd05} for 
$d\ge 3$,  is the dimensional threshold $s_0=\frac{d}{2}+\frac{1}{3}$. If the set $E$ is assumed to be Ahlfors-David regular, then the 
dimensional threshold $s_0=1$ was recently nearly established  for $d=2$ by Orponen \cite{O15}. More precisely, he showed that if the 
Hausdorff dimension of an Ahfors-David regular set in the plane is $1$, then the upper Minkowski dimension of the distance set is also $1$. An 
example due to Falconer \cite{Fal86} shows  that $s_0=1$ is essentially sharp:  if the Hausdorff dimension of a planar set is less than one, then 
the upper Minkowski dimension of the distance set can in general be less than 1. 

\subsection{Congruence and similarity classes of simplexes} The distance problem can be viewed as a question about two-point 
configurations. A pair of points $x,y \in E$ is congruent to another pair $x',y' \in E$ iff $|x-y|=|x'-y'|$. This induces an equivalence relation $\sim$ 
and we may view $\Delta(E)$ as $E \times E \ \backslash \sim$. This set can be naturally identified with the distance set $\Delta(E)$. In the 
same way, we may consider $(k+1)$-point configurations $\{x^1, x^2, \dots, x^{k+1}\}$, $x^j \in E$, $k \leq d$, 
and  say that $\{x^1, \dots, x^{k+1}\}$ is {\em congruent} to $\{y^1, \dots, y^{k+1}\}$ if there exists a translation $\tau \in {\Bbb R}^d$ and a 
rotation $g \in O_d({\Bbb R})$ such that $y^j=\tau+gx^j, \ j=1,2, \dots, k+1$. 
The resulting equivalence relation allows us to define $T_k(E):=E \times E \times \dots \times E \ 
\backslash \sim$, which can be  identified with a ${k+1 \choose 2}$-tuple of distances $|x^i-x^j|$, $1 \leq i<j \leq k+1$. We can also consider 
$T_k(E)$ as the  set of equivalence classes of $k$-dimensional simplexes determined by points of $E$, and can be viewed as a subset of ${\Bbb R}^{k+1 \choose 2}$ for the purpose of measuring its size in terms of Hausdorff dimension of non-vanishing Lebesgue measure.\footnote{For these considerations, one can ignore the action of the permutation group on $(x^1,\dots,x^{k+1})$, since this does not affect the non-vanishing of the ${k+1 \choose 2}$-dimensional Lebesgue measure.}

The natural generalization of the Falconer distance problem in this context is to find a dimensional threshold $s_0$ such that if the  Hausdorff 
dimension, $\hdim(E)$, of a compact set $E$ is greater than $s_0=s_0(k,d)$, then the ${k+1 \choose 2}$-dimensional Lebesgue measure of 
$T_k(E)$ is positive: ${\mathcal L}^{k+1 \choose 2}(T_k(E))>0$. A variety of results have been obtained in this direction in recent years using 
everything from multi-linear theory to group actions. See, for example, \cite{GI12,EHI12,GGIP12}. To various extents those papers were 
preceded and motivated by finite field models worked out in \cite{BHIPR13,BIP12,CEHIK12, HI07}. As a result of these efforts we know that for a 
compact set $E \subset {\Bbb R}^d$ eith  $d \ge 3$ and $2 \leq k \leq d$,  if  
\begin{equation} \label{ehiexp}
\hdim(E)> \min \left\{\frac{dk+1}{k+1}, \frac{d+k+1}{2} \right\},
\end{equation}
then  ${\mathcal L}^{k+1 \choose 2}(T_k(E))>0$. 
The purpose of the current paper is to improve (i.e., lower) the $\frac{d+k+1}{2}$ threshold obtained in \cite{EHI12} down to $\frac{d+k}{2}$, and introduce a new method in doing so.
\medskip

Closely related to these questions are estimates for multilinear forms which, borrowing a term for analogous expressions in the discrete setting,
we call \emph{incidence bounds}. 
For $\rho\in C^\infty_0(\R^d)$ with $\rho\ge 0$, 
$supp(\rho)\subset \{|x|\le 1\}$, $\rho\equiv 1$ on  $\{|x|\le\frac{1}{4}\}$ and $\| \rho \|_{L^1(\mathbb{R}^d)} = 1$,
form the approximate identity $\rho^{\epsilon}(x)=\epsilon^{-d}\rho(x/\epsilon)$.
For $t>0$, let $\sigma_t$ denote the surface measure on the sphere of radius $t$,
and $\sigma_t^\epsilon=\sigma_t*\rho_\epsilon$. Finally, let  $\mu$ be a Frostman measure supported on a set $E \subset {\Bbb R}^d$, $d \ge 2$ (see \cite{Mat95}).
For positive numbers $\{t_{ij}\}$, we ask whether an incidence bound,
\begin{equation} \label{incidencebound} 
|\Lambda^{\epsilon}_{k,d}(\mu)|:=\Big|\int \dots \int \prod_{1 \leq i<j \leq k+1} \sigma_{t_{ij}}^{\epsilon}(x^i-x^j) \prod_{l=1}^{k+1} d\mu(x^l)\Big|\lesssim 1,  \end{equation} 
holds.
(Here and throughout, $X \lesssim Y$ means that there exists $C>0$ such that $X \leq Y$ independent of $\epsilon$.) 

Whenever (\ref{incidencebound} ) holds, 
it  implies  that the ${k+1 \choose 2}$-dimensional Lebesgue measure of $T_k(E)$ (defined above) is positive (see \cite{GGIP12}). But the
uniform estimate in (\ref{incidencebound}) is considerably stronger than the positivity of the Lebesgue measure. For example, in  \cite{GI12}, the authors proved that, in the case $k=d=2$, 
(\ref{incidencebound}) holds  if  $\hdim(E)>\frac{7}{4}$, yielding not only the continuous Falconer-type configuration result but also a discrete  
result:  If $A \subset {\Bbb R}^2$ is a finite 
homogeneous set  with $|A|=N$, then the number of triples of points from $A$ determining an equilateral triangle of fixed side length does not exceed 
$CN^{\frac{9}{7}}$,  an improvement over the previously known $Cn^{\frac{4}{3}}$ bound (which is a consequence of the Szemeredi-Trotter incidence theorem). (For applications of continuous incidence bounds in geometric measure theory, see, e.g., \cite{EIT13}; for the definition of homogeneous set, see  \cite{SoVu}.) 

\vskip.25in 

\subsection{Statement of results} 

The main results of this paper are the following. 

\begin{theorem} \label{main} 
Let $E$ be a compact subset of ${\Bbb R}^d$, $d \ge 3$, $2 \leq k \leq d$. Suppose that  
\begin{equation} \label{genbest} \hdim(E)> \min \left\{\frac{dk+1}{k+1}, \frac{d+k}{2} \right\}. \end{equation} 
Then ${\mathcal L}^{k+1 \choose 2}(T_k(E))>0$. 
\end{theorem} 

\vskip.125in 

\begin{remark} The following observations clarify the role of the exponents in Theorem \ref{main}. 

\begin{itemize} 

\item The estimate (\ref{genbest}) improves(i.e., lowers) by $\frac{1}{2}$ the best exponent previously known (proved in \cite{EHI12} and described in (\ref{ehiexp})). 

\item The sufficiency as a lower bound of the first term in the $\min$  was previously obtained in \cite{GILP13}.

\item For $k=d=2$, the best lower bound known is $\hdim(E)>\frac{8}{5}$. 

\end{itemize} \end{remark}

\vskip.125in 

\begin{theorem} (Incidence Bound) \label{main2} Let $E$ be a compact subset of ${\Bbb R}^d$, $d \ge 2$, $2 \leq k \leq d$ and let $\mu$ be a Frostman measure on $E$. Then (\ref{incidencebound}) holds if 
$$ dim_{{\mathcal H}}(E)> \min \left\{ \frac{d+k}{2}, d-\frac{d-1}{2k} \right\}.$$  
\end{theorem} 

\vskip.125in 

\begin{remark} \label{littleleft} The  sufficiency as a lower bound of the second term in the $\min$ was previously obtained in \cite{GGIP12}. The case $d=k=2$ was handled earlier in \cite{GI12}. 
Note that $d-\frac{d-1}{2k}<\frac{d+k}{2}$ if and only if $k=d$. \end{remark} 

\vskip.125in 

\begin{remark} Let $\alpha_{k,d}$ denote the optimal exponent for the congruent $d$-dimensional simplex problem, i.e., $\alpha_{k,d}$ is the infimum of those $\alpha$ for which $\mathcal{L}^{k+1 \choose 2}(T_k(E)) > 0$ whenever $\hdim(E)>\alpha$. Easy examples show that $\alpha_{k,d}\geq \max\left\lbrace k-1, \frac{d}{2} \right\rbrace$. For non-trivial sharpness examples Burak Erdo\u{g}an and the second listed author obtained $\alpha_{2,2}\geq \frac{3}{2}$, 
while Jonathan DeWitt, Kevin Ford, Eli Goldstein, Steven J. Miller, Gwyneth Moreland, the fourth listed author and Steven 
Senger obtained $\alpha_{k,d}\geq\frac{d(k+1)}{d+2}$ which is only non-trivial in the range $\frac{d}{2} < k\leq d$. The first 
novel case in Theorem \ref{main}, where the new threshold $\frac{d+k}{2}$ beats the previously obtained threshold $\frac{dk+1}{k+1}$, is when $d=5$, $k=2$. There we obtain the threshold $\frac{7}{2}$, which we can only compare to the trivial sharpness 
example $\alpha_{2,5} \geq 3$. The first new case where $\frac{d+k}{2}$ beats out the old threshold and where we have $k>\frac{d}{2}$ (in order to have a non-trivial sharpness example) is for $d=7$, $k=4$. In this case Theorem \ref{main} yields the threshold $\frac{11}{2}$, while the non-trivial sharpness example shows $\alpha_{4,7} \geq \frac{35}{9}$. A gap remains. \end{remark}

\vskip.125in 

The following multi-linear estimate follows easily from the proof of Theorem \ref{main2}. 

\begin{theorem} For $\epsilon>0$, let $\sigma^\epsilon=\sigma_1*\rho_\epsilon$, and define the multi-linear operator
$$ M^{\epsilon}_k(f_1, \dots, f_k)(x^{k+1})=\int \dots \int \prod_{1 \leq i<j \leq k+1} \sigma^{\epsilon}(x^i-x^j) \prod_{l=1}^{k} 
f_l(x^l)\,d\mu(x^l).$$
Then if $\mu$ is a compactly supported Frostman measure on a set $E\subset\R^d$ with
$\hdim(E)>\min \left\{\frac{d+k}{2}, d-\frac{d-1}{2k}\right\}$, then 
$$ M^{\epsilon}_k: L^k(\mu) \times \dots \times L^k(\mu) \to L^1(\mu),$$ 
with constants independent of $\epsilon$. 
\end{theorem}

\section{Proof of Theorem \ref{main}} 

\vskip.125in

As mentioned in the introduction, the fact that ${\mathcal L}^{k+1 \choose 2}(T_k(E))>0$ if  $\hdim(E)>\frac{dk+1}{k+1}$ was established in \cite{GILP13}. We are thus left with establishing the $\frac{d+k}{2}$ threshold. 

Begin with the following observation proved in \cite{GGIP12}. Heuristically, it says that if a given a configuration does not arise too often, there must be many different configurations. In the language of combinatorics, this is an incidence estimate. 

\begin{lemma} \label{reductionlemma} 
Let $E\subset\R^d,\, d\ge 2,\, 2\le k\le d$, $\mu$ be a Frostman measure on $E$, and $\sigma_t,\, \sigma_t^\epsilon$ be as above. Let $t_{ij}>0$ be arbitrary.
Suppose that 
\begin{equation} \label{mamaform} \big|\Lambda_k(\mu)\big|:=\Big|\int \dots \int \prod_{1 \leq i<j \leq k+1} \sigma^{\epsilon}_{t_{ij}}(x^i-x^j) d\mu(x^1) \dots d\mu(x^{k+1})\Big| \lesssim 1 \end{equation} 
with constants independent of $\epsilon$. 
Then ${\mathcal L}^{k+1 \choose 2}(T_k(E))>0$. 
\end{lemma} 

\begin{remark} Observe that $\Lambda_k(\mu)$ depends on $\epsilon$ and we shall obtain bounds independent of $\epsilon$. \end{remark}

\begin{remark} This incidence estimate also establishes the first threshold in Theorem \ref{main2}. \end{remark} 

\subsection{The case $k=2$} We now bound (\ref{mamaform}) in the case $k=2$ because it is slightly simpler and illustrates the method. For the sake of notational simplicity, we write out the case $t_{12}=t_{13}=t_{23}=1$, 
but the argument works  in the general case of general $t_{ij}$. 
Denote $\sigma_1^\epsilon$ by $\sigma^\epsilon$.
We have 
$$ \Lambda_2(\mu)=\int \int \int \sigma^{\epsilon}(x-y)\sigma^{\epsilon}(x-z)\sigma^{\epsilon}(y-z)d\mu(x)d\mu(y)d\mu(z).$$ 

Given complex numbers $\alpha,\beta,\gamma$, define 
$$ \Lambda_2^{\alpha,\beta,\gamma}(\mu)=\int \int \int \sigma^{\epsilon,\alpha}(x-y) \sigma^{\epsilon,\beta}(x-z)\sigma^{\epsilon,\gamma}(y-z) d\mu(x)d\mu(y)d\mu(z),$$
where
\begin{equation*}
\sigma^{\epsilon,z}(x) := \frac{2^{\frac{d-z}{2}}}{\Gamma\left(z/2 \right)}(\sigma^{\epsilon} * |\cdot |^{-d+z})(x)
\end{equation*}
is initially defined for $\text{Re}(z) > 0$ but then extended to the complex plane by analytic continuation. This follows the strategy introduced by the first two authors in \cite{GI12}.

\vskip.125in 

\begin{theorem} \label{2form} Let $E$, $\mu$ and $\Lambda_2^{\alpha,\beta,\gamma}$ be as above. 

\vskip.125in 

Suppose that $Re(\alpha)=1+\delta$, $Re(\beta)=Re(\gamma)=-\frac{1}{2}-\frac{\delta}{2}$ for some small $\delta>0$. Then 
\begin{equation} \label{lambdaest} |\Lambda_2^{\alpha,\beta,\gamma}(\mu)| \lesssim 1 \end{equation} 
with constants independent of $\epsilon$ provided that $dim_{{\mathcal H}}(E)>\frac{d+2}{2}$. 

\end{theorem} 

\begin{corollary} \label{3lineslemmaapp} In view of Theorem \ref{2form} and the three lines lemma, (\ref{lambdaest}) holds if $\alpha=\beta=\gamma=0$ and $dim_{{\mathcal H}}(E)>\frac{d+2}{2}$. \end{corollary} 

\vskip.125in 

To prove (\ref{lambdaest}) we shall need the following result due to the second named author, Krause, Sawyer, Taylor and Uriarte-Tuero \cite{IKSTU}.  

\begin{theorem} \label{maintool} Let $\mu, \nu$ be compactly supported Borel measures satisfying $\mu(B(x,r)) \leq Cr^{s_{\mu}}$, $\nu(B(x,r)) \leq Cr^{s_{\nu}}$, and $\lambda$ compactly supported Borel measure satisfying $|\widehat{\lambda}(\xi)| \leq C{|\xi|}^{-\eta}$. Define $T_{\lambda^{\epsilon}}f(x)=\int \lambda^{\epsilon}(x-y) f(y) d\mu(y)$. Let $s=\frac{s_{\mu}+s_{\nu}}{2}$, and suppose that $\eta>d-s$. Then 
\begin{equation} \label{ytm} {||T_{\lambda^{\epsilon}}f||}_{L^2(\nu)}  \leq C{||f||}_{L^2(\mu)} \end{equation} with constant $C$ independent of $\epsilon$.  
\end{theorem}

We include the proof of Theorem \ref{maintool}, for the sake of completeness, in Section \ref{sectionmaintool} below. With Theorem \ref{maintool} in tow, the proof of Theorem \ref{2form} proceeds as follows. Let $Re(\alpha)=1+\frac{\delta}{2}$. Then
\begin{equation}\label{AfterFirstStep}
|\Lambda_2^{\alpha,\beta,\gamma}(\mu)| \leq C \int \int \int |\sigma^{\epsilon,\beta}|(x-z) |\sigma^{\epsilon, \gamma}|(y-z) 
d\mu(x) d\mu(y) d\mu(z).
\end{equation}

Here we need the following basic calculation. 

\begin{lemma}\label{Linfty}
Let $\sigma^{\epsilon,\alpha}$ be as above. Then
$$ |\sigma^{\epsilon,\alpha}(x-y)| \lesssim 1 .$$
\end{lemma}
To prove this lemma first write
$$ \sigma^{\epsilon,\alpha}(x-y) = \frac{2^{\frac{d-\alpha}{2}}}{\Gamma\left(\alpha/2 \right)}\int \sigma^{\epsilon}(\tau)|x-y-\tau|^{-d+\alpha}d\tau,$$
and then decompose the integral in $\tau$ into annuli where $|x-y-\tau|\sim 2^{-j}$ and note that since we are in a compact setting then this happens for a set of $j$ that are bounded below, say $-M\leq j$, where $M$ only depends on the diameter of the set $E$. Finally conclude
\begin{align*}
|\sigma^{\epsilon,\alpha}(x-y)| &\lesssim \frac{2^{\frac{d-Re(\alpha)}{2}}}{|\Gamma\left(\alpha/2 \right)|}\sum\limits_{j=-M}^{\infty}\ 2^{j(d-Re(\alpha))} \int\limits_{|x-y-\tau| \sim 2^{-j}}\sigma^{\epsilon}(\tau)d\tau \\
&\lesssim \frac{2^{\frac{d-Re(\alpha)}{2}}}{|\Gamma\left(\alpha/2 \right)|}\sum\limits_{j=-M}^{\infty} 2^{j(d-Re(\alpha))} 2^{-j(d-1)} \\
&\lesssim \frac{2^{\frac{d-Re(\alpha)}{2}}}{|\Gamma\left(\alpha/2 \right)|}\sum\limits_{j=-M}^{\infty} \left( 2^{\frac{\delta}{2}} \right)^{-j} \lesssim 1
\end{align*}
This completes the proof of Lemma \ref{Linfty}.

Applying Cauchy-Schwarz to (\ref{AfterFirstStep}) reduces matters to bounding 
\begin{equation} \label{aftercs} \int {||\sigma^{\epsilon, \beta}|*\mu(x)|}^2 d\mu(x) \end{equation} and the same expression with $\gamma$ instead of $\beta$. It can be deduced from the classical stationary phase estimates (see, e.g., \cite{St93}) that 
$$ |\widehat{|\sigma^{\epsilon, \beta}|}(\xi)| \leq C{|\xi|}^{-\frac{d-1}{2}-Re(\beta)}.$$

In view of this and the assumption that $\mu$ in the definition of $\Lambda^{\alpha,\beta,\gamma}_2(\mu)$ is a Frostman measure, we may apply Theorem \ref{maintool} with $\mu=\nu$, $s_{\mu}$ slightly smaller than  $\hdim(E)$ and $\eta=\frac{d-1}{2}-Re(\beta)=\frac{d-2}{2}-\frac{\delta}{2}$. It follows that the expression in (\ref{aftercs}) is bounded if $\hdim(E)>d-\frac{d-2}{2}$, as claimed. 

\vskip.125in 

\subsection{The general case} In the case $k=2$ we reduced matters to the chain of length $2$. In general, we shall reduce matters to chains of length $k$. Let 
$$ \Lambda_k^{\alpha}(\mu)=\int \dots \int \prod_{1 \leq i<j \leq k+1} \sigma^{\epsilon,\alpha_{ij}}_{t_{ij}}(x^i-x^j) d\mu(x^1) \dots d\mu(x^{k+1}).$$

Now set ${k \choose 2}$ of the $Re(\alpha_{ij})$s equal to $1+\delta {k \choose 2}^{-1}$, where the $\alpha_{ij}$s are chosen in such a way that what remains is a $k$-link chain. Choose the remaining $\alpha_{ij}$s equal to $-\frac{k-1}{2}-\frac{\delta}{k}$. This is accomplished as follows. There are $(k+1)!$ ways to order  $1,2, \dots, k+1$. Each such ordering can be written as
$(N(1),N(2), \dots, N(k+1))$, where $N$ is a bijection on $\{1,2, \dots, k+1\}$. Given such an ordering, we set $Re(\alpha_{ij})=-\frac{k-1}{2}-\frac{\delta}{k}$ if $(i,j)=(N(m), N(m+1))$ or 
$(N(m+1),N(m))$ (depending on whether or not $N(m)<N(m+1)$). For the remaining $(i,j)$, $1 \leq i<j \leq k+1$, set $Re(\alpha_{ij})=1+\delta {k \choose 2}^{-1}$. 

It follows from Lemma \ref{Linfty} that, up to a relabeling of vertices, 
\begin{equation} \label{zeelow} |\Lambda_k^{\beta}(\mu)| \leq \int \dots \int \prod_{j=1}^k |\sigma_{t_j}^{\epsilon, \beta_j}|(x^{j+1}-x^j) 
d\mu(x^1) \dots d\mu(x^{k+1}),\end{equation} where 
$$ Re(\beta_j)=-\frac{k-1}{2}-\frac{\delta}{k}.$$ 

We shall see using a modification of an argument in \cite{BIT14} that the right hand side of (\ref{zeelow}) is bounded if  $\hdim(E)>\frac{d+k}{2}$. We shall need the following generalization of an upper bound from \cite{BIT14}. 

\begin{theorem} \label{chaingen} Let $\lambda_1,\ldots,\lambda_n$ be compactly supported Borel measures such that 
$$ |\widehat{\lambda_j}(\xi)| \leq C{|\xi|}^{-\alpha} \ \text{for some} \ \alpha>0 \text{ and all } 1\leq j \leq n.$$ 

Let $\lambda_j^{\epsilon}(x)=\lambda_j*\rho_{\epsilon}(x)$ where $\rho_{\epsilon}$ is as above. Let $\mu$ be a Frostman measure on a compact set of Hausdorff dimension $s$. Then 
\begin{equation}\label{lambdaineq}
\left| \int \dots \int \prod_{j=1}^n \lambda_j^{\epsilon}(x^{j+1}-x^j) d\mu(x^1) \dots d\mu(x^{n+1}) \right| \leq C^{n+1}
\end{equation}
independent of $\epsilon$ if $s>d-\alpha$. Here $C$ is the constant obtained in Theorem \ref{maintool}.
\end{theorem}

Using the fact that
$$ |\widehat{\sigma^{\epsilon,\beta_j}_{t_j}}(\xi)| \lesssim |\xi|^{-\frac{d-1}{2}-Re(\beta_j)} $$
where $Re(\beta_j)=-\frac{k-1}{2}-\frac{\delta}{k}$ allows us to conclude,
using Theorem \ref{chaingen}, that the right hand side of (\ref{zeelow}) is bounded if
$$ \hdim E > d - \left(\frac{d-1}{2} - \frac{k-1}{2}-\frac{\delta}{k}\right) = \frac{d+k}{2} + \frac{\delta}{k} $$
for any $\delta > 0$. Finally observe that the sum of the $Re(\alpha_{ij})$ is $0$ which shows that $\vec{0}$ is in the convex hull of all the points obtained by permuting the $\alpha_{ij}$ and this completes the proof of Lemma \ref{reductionlemma} via the three lines lemma.

To prove Theorem \ref{chaingen} we follow a similar strategy as in \cite{BIT14}. Let
$$ T_k^{\epsilon}f(x) = \left(\lambda_k * (f\mu) \right)(x), $$
$f_0^{\epsilon}(x) = 1$, and then recursively define $ f_{k+1}^{\epsilon} = T^{\epsilon}_{k+1}f^{\epsilon}_k,\, k\ge 0$.
With this notation we can write (\ref{lambdaineq}) as
$$ \left|\int f_{n+1}^{\epsilon}(x^{n+1}) d\mu(x^{n+1})\right| \leq C .$$
We start with the left hand side and apply the Cauchy-Schwarz inequality,
$$ \left|\int f_{n+1}^{\epsilon}(x^{n+1}) d\mu(x^{n+1})\right| \leq \| f_{n+1}^{\epsilon} \|_{L^2(\mu)}, $$
where we use that $\mu$ is a probability measure so $\int d\mu(x^{n+1}) = 1$. Then the proof is completed with repeated use of Theorem \ref{maintool}:
$$ \| f_{n+1}^{\epsilon} \|_{L^2(\mu)} \leq C \| f_{n}^{\epsilon} \|_{L^2(\mu)} \leq \ldots \leq C^{n+1} \| f_{0}^{\epsilon} \|_{L^2(\mu)} = C^{n+1} $$
where $\mu = \nu$, with $s_{\mu}$ slightly smaller than $\hdim(E)$ and $\alpha > d-s$.

\vskip.25in 

\section{Proof of Theorem \ref{maintool}} 
\label{sectionmaintool}

It is enough to show that if $g \in L^2(\nu)$, then 
\begin{equation} \label{setup} |<T_{\lambda^{\epsilon}}f, g\nu>| \ \leq C{||f||}_{L^2(\mu)} \cdot {||g||}_{L^2(\nu)}, \end{equation} where the constant $C$ is independent of $\epsilon$. 

\vskip.125in 

The left hand side of (\ref{setup}) equals 
\begin{equation} \label{preinversion} \int \lambda^{\epsilon}*(f\mu)(x) g(x) d\nu(x). \end{equation} 

Indeed, 
$$ \lambda^{\epsilon}*(f\mu)(x)=\int e^{2 \pi i x \cdot \xi} \widehat{\lambda}(\xi) \widehat{\rho}(\epsilon \xi)
\widehat{f\mu}(\xi) d\xi$$ for every $x \in {\Bbb R}^d$ because the left hand side is a continuous $L^2({\Bbb R}^d)$ function and $$\widehat{\lambda}(\cdot) \widehat{\rho}(\epsilon \cdot) \widehat{f\mu}(\cdot) \in L^1 \cap L^2({\Bbb R}^d).$$

It follows that (\ref{preinversion}) equals 

$$\int \int e^{2 \pi i x \cdot \xi} \widehat{\lambda}(\xi) \widehat{\rho}(\epsilon \xi) \widehat{f \mu}(\xi) d\xi g(x) d\nu(x).$$ Applying Fubini, we see that this expression equals 
$$ \int \int e^{2 \pi i x \cdot \xi} g(x) d\nu(x) \widehat{\lambda}(\xi) \widehat{\rho}(\epsilon \xi) \widehat{f \mu}(\xi) d\xi$$
$$=\int \widehat{\lambda}(\xi) \widehat{\rho}(\epsilon \xi) \widehat{f \mu}(\xi)  \widehat{g \nu}(\xi) d\xi.$$
The modulus of this expression is bounded by an $\epsilon$-independent multiple of 

$$ \int {|\xi|}^{-\alpha} |\widehat{f\mu}(\xi)| \cdot  |\widehat{g \nu}(\xi)| d\xi.$$ 

By Cauchy-Schwarz, this expression is bounded by 

\begin{equation} \label{almostthere} {\left( \int {|\widehat{f\mu}(\xi)|}^2 {|\xi|}^{-\alpha_{\mu}} d\xi \right)}^{\frac{1}{2}} \cdot {\left( \int {|\widehat{g\nu}(\xi)|}^2 {|\xi|}^{-\alpha_{\nu}} d\xi \right)}^{\frac{1}{2}}=:\sqrt{A} \cdot \sqrt{B}, \end{equation} where $\alpha_{\mu}, \alpha_{\nu}>0$ and $\frac{\alpha_{\mu}+\alpha_{\mu}}{2}=\alpha$. 

\vskip.125in 

\begin{lemma} \label{Aestimation} With the notation above, we have 

$$ A \leq C{||f||}^2_{L^2(\mu)}\hbox{ and }B \leq C{||g||}^2_{L^2(\nu)}$$ if 
\begin{equation} \label{numuconditions}\alpha_{\mu}>d-s_{\mu}\hbox{ and } \alpha_{\nu}>d-s_{\nu}. 
\end{equation} 
\end{lemma} 

\vskip.125in 

We give a direct argument in the style of the proof of Lemma 7.4 in \cite{W04}. It is enough to prove that the estimate for $A$ follows from the condition on $\alpha_{\mu},\, s_\mu$, since the estimate for $B$ follows from the same statement applied to $\alpha_\nu,\, s_\nu$. By Proposition 8.5 in \cite{W04}, 
\begin{equation} \label{schursetup} A=\int \int f(x)f(y) {|x-y|}^{-d+\alpha_{\mu}} d\mu(x)d\mu(y)=<f,Uf>, \end{equation} where 

$$ Uf(x)=\int {|x-y|}^{-d+\alpha_{\mu}} f(y) d\mu(y).$$ 

Observe that 
$$ \int {|x-y|}^{-d+\alpha_{\mu}} d\mu(y)= \int {|x-y|}^{-d+\alpha_{\mu}} d\mu(x)$$
$$ \leq C \sum_{j>0} 2^{j(d-\alpha_{\mu})} \int_{2^{-j} \leq |x-y| \leq 2^{-j+1}} d\mu(y)$$
$$ \leq C' \sum_{j>0} 2^{j(d-\alpha_{\mu}-s_{\mu})} \leq C'' \ \text{if} \ \alpha_{\mu}>d-s_{\mu}.$$

It follows by Schur's test (see Lemma 7.5 in \cite{W04} and the original argument in \cite{Schur11}) that 
$${||Uf||}_{L^2(\mu)} \leq C'' {||f||}_{L^2(\mu)}$$ and we are done in view of (\ref{schursetup}) and Cauchy-Schwarz. 

\vskip.25in 

\section{Proof of the second estimate in Theorem \ref{main2}} 
\label{bookproof} 

\vskip.125in

As mentioned in the introduction, the fact that ${\mathcal L}^{k+1 \choose 2}(T_k(E))>0$ if  $\hdim(E)>d-\frac{d-1}{2k}$ was established in \cite{GGIP12}. We give a simpler, more transparent proof below.

As we noted in Remark \ref{littleleft}, we only need to deal with the case $k=d$, since
for $k\le d-1$, the known $\alpha_{k,d}\le \frac{d+k}2$ is the minimum of the two quantities in Theorem \ref{main2}. 
We prove the estimate in the case $t_{ij}=1$ as the proof of the general case is essentially the same. Let $\mu$ be a Frostman measure on $E$ and define 
$\alpha=(\alpha_1, \alpha_2, \dots, \alpha_{d+1})$, where $\alpha_j \in {\mathbb C}$. Define $\Lambda^{\delta}(\alpha)$ by the integral
$$ \int \int \sigma^{\epsilon}(x^1-x^2) d\mu^{\delta}_{\alpha_1}(x^1) d\mu^{\delta}_{\alpha_2}(x^2) \int \dots \int 
\prod_{1 \leq i<j \leq d+1; (i,j) \not=(1,2)} \sigma^{\epsilon}(x^i-x^j) \prod_{l=3}^{d+1} d\mu^{\delta}_{\alpha_l}(x^l),$$ where 
$$\mu^{\delta}(x)=\mu*\rho_{\delta}(x),$$ with $\rho$ a smooth cut-off function, $\int \rho=1$, $\rho_{\delta}(x)=\delta^{-d} \rho(x/\delta)$, and 
\begin{equation} \label{zkaif} \mu_z^{\delta}(x) := \frac{2^{\frac{d-z}{2}}}{\Gamma\left(z/2 \right)}(\mu^{\delta} * |\cdot |^{-d+z})(x). \end{equation}

We shall prove that when $Re(\alpha_j)$, $j=3, \dots, d+1$, equals $\frac{d-1}{2d}-\frac{\epsilon}{d-1}$, then 
\begin{equation} \label{intersectionstep} \left| \int \dots \int \prod_{1 \leq i<j \leq d+1; (i,j) \not=(1,2)} \sigma^{\epsilon}(x^i-x^j) 
\prod_{l=3}^{d+1} d\mu^{\delta}_{\alpha_l}(x^l) \right| \leq C \end{equation} provided that 
$\hdim(E)>d-\frac{d-1}{2d}+\frac{\epsilon}{d-1}$. This will follow from two observations. First, we are going to show that 
\begin{equation} \label{boundeddude} |\mu^{\delta}_z(x)| \leq C \end{equation} provided that $dim_{{\mathcal H}}(E)>d-Re(z)$. 

We will then show that 
\begin{equation} \label{elemgeometry} \left| \int \dots \int 
\prod_{1 \leq i<j \leq d+1; (i,j) \not=(1,2)} \sigma^{\epsilon}(x^i-x^j) dx_3 \dots dx_{d+1} \right| \leq C \end{equation} using elementary geometric considerations. Combining (\ref{elemgeometry}) and (\ref{boundeddude}) yields (\ref{intersectionstep}). 

\vskip.125in 

We shall then prove that if $Re(\alpha_1)=Re(\alpha_2)=-\frac{{(d-1)}^2}{4d}+\frac{\epsilon}{2}$, then 
\begin{equation} \label{falconer+} \int \int \sigma^{\epsilon}(x^1-x^2) d|\mu_{\alpha_1}|(x^1)d|\mu_{\alpha_2}|(x^2) \leq C \end{equation} under the same constraint on $dim_{{\mathcal H}}(E)$ as above. 

\vskip.125in 

Combining (\ref{falconer+}) and (\ref{elemgeometry}) will show that 
$$ |\Lambda^{\delta}(\alpha)| \leq C$$ if $Re(\alpha_1)=Re(\alpha_2)=-\frac{{(d-1)}^2}{4d}+\frac{\epsilon}{2}$ and $Re(\alpha_j)=\frac{d-1}{2d}-\frac{\epsilon}{d-1}$, $j=3, \dots, d+1$. Interchanging the role of variables, we obtain ${d+1 \choose 2}$ relations of this type and since the real parts add up to $0$ we have $\vec{0}$ in the convex hull of all these relations and thus the result easily follows by analytic interpolation. 

We now estblishes estimates (\ref{boundeddude}), (\ref{elemgeometry}), (\ref{falconer+}). 

\subsection{Proof of the estimate (\ref{boundeddude})} 

We have 
\begin{align*}
|\mu_z^{\delta}(x)| &=\left| \frac{2^{\frac{d-z}{2}}}{\Gamma\left(z/2 \right)} \right| \int {|x-y|}^{-d+z} d\mu_{\delta}(y) \\
&=\left| \frac{2^{\frac{d-z}{2}}}{\Gamma\left(z/2 \right)} \right| \sum_j \int_{2^{-j-1}<|x-y| \leq 2^{-j}} {|x-y|}^{-d+z} d\mu^{\delta}(y) \\
& \leq C \left| \frac{2^{\frac{d-z}{2}}}{\Gamma\left(z/2 \right)} \right| \sum_j 2^{j(d-Re(z))} \int_{2^{-j-1}<|x-y| \leq 2^{-j}} d\mu^{\delta}(y)\\
& \leq C_{\alpha}'(z) \sum_j 2^{j(d-Re(z)-\alpha)}
\end{align*}
for any $\alpha<dim_{{\mathcal H}}(E)$ since $\mu$ is a Frostman measure. The geometric series converges if $Re(z)>d-dim_{{\mathcal H}}(E)$, as claimed. This completes the proof of the estimate (\ref{boundeddude}). 

\vskip.125in 

\subsection{Proof of the estimate (\ref{elemgeometry})} This is a volume packing estimate. We have a transverse intersection of $\epsilon$ neighborhoods of ${d+1 \choose 2}-1$ surfaces in ${\Bbb R}^{d(d-2)}$ multiplied by $\epsilon^{-{d+1 \choose 2}+1}$. This results in an $\epsilon$ neighborhood of an $d(d-2)-{d+1 \choose 2}+1$ dimensional surface, multiplied by  $\epsilon^{-{d+1 \choose 2}+1}$. The resulting weighted volume is $\approx 1$ and the proof is complete. 

\vskip.125in 

\subsection{Proof of the estimate (\ref{falconer+})} This result follows from the following general observation combined with (\ref{zkaif}). 
\begin{lemma} Let $d\mu_A, d\mu_B$ denote compactly supported Borel measures and $s_A,s_B$ real numbers such such that 
$$ \max \left\{ \int {|\widehat{\mu}_A(\xi)|}^2 {|\xi|}^{-d+s_A} d\xi, \int {|\widehat{\mu}_B(\xi)|}^2 {|\xi|}^{-d+s_B} d\xi  \right\} \leq C.$$ 

Let $\lambda \in L^1({\Bbb R}^d)$ such that 
$$ |\widehat{\lambda}(\xi)| \leq C{|\xi|}^{-\gamma} \ \text{for some} \ \gamma>0.$$ 

Then 
$$ \int \int \lambda(x-y) d\mu_A(x)d\mu_B(y) \leq C' \ \text{if} \ \frac{s_A+s_B}{2}>d-\alpha.$$ 
\end{lemma} 

\vskip.125in 

The proof is very similar to that of Theorem \ref{maintool}. Choose $\alpha_A, \alpha_B$ such that $\frac{\alpha_A+\alpha_B}{2}=\alpha$ and write 
$$  \int \int \lambda(x-y) d\mu_A(x)d\mu_B(y)=\int \widehat{\mu}_A(\xi) \overline{\widehat{\mu}_B(\xi)} \widehat{\lambda}(\xi) d\xi$$ 
$$ \leq \int |\widehat{\mu}_A(\xi)|  |\widehat{\mu}_B(\xi)| {|\xi|}^{-\alpha} d\xi=\int |\widehat{\mu}_A(\xi)| {|\xi|}^{-\frac{\alpha_A}{2}}  
|\widehat{\mu}_B(\xi)| {|\xi|}^{-\frac{\alpha_B}{2}} d\xi$$ 
$$ \leq {\left( \int {|\widehat{\mu}_A(\xi)|}^2 {|\xi|}^{-\alpha_A} d\xi \right)}^{\frac{1}{2}} \cdot {\left( \int {|\widehat{\mu}_B(\xi)|}^2 {|\xi|}^{-\alpha_B} 
d\xi \right)}^{\frac{1}{2}}=I \cdot II.$$

By assumption, $I$ is bounded if $\alpha_A \ge d-s_A$ and $II$ is bounded if $\alpha_B \ge d-s_B$. This completes the proof of the lemma and hence the proof of the second estimate in Theorem \ref{main2}. 


\bigskip
\bigskip

\end{document}